\documentclass[12pt]{article} 
\usepackage{amsmath,amsthm, amssymb} 
\usepackage{amssymb,latexsym}

\newtheorem{theorem}{Theorem}

\newtheorem{lemma}{Lemma}

\begin{document}

\baselineskip=17pt 

\title{\bf The ternary Goldbach problem \\ with arithmetic weights attached \\ to two of the variables}

\author{\bf D. I. Tolev \footnote{Supported by Sofia University Grant 221}}

\date{}
\maketitle

\begin{abstract}
We consider the ternary Goldbach problem with two prime variables of the form $k^2+m^2+1$
and find an asymptotic formula for the number of its solutions.
\end{abstract}

\paragraph{Notations.}

By greek letters we denote real numbers and 
by small latin letters --- integers. However, the letter $p$, with or without subscripts, is reserved for primes.
By $\varepsilon$ we denote an arbitrarily small positive number, not the same in different appearances.
$N$ is a sufficiently large odd integer and $\mathcal L = \log N$. 
We denote by $\mathcal J$
the set of all subintervals of the interval $[1, N]$ and
if $J_1, J_2 \in \mathcal J$ then 
${\bf J} = \langle J_1, J_2\rangle $ is the corresponding
ordered pair. Respectively
${\bf k} = \langle k_1, k_2 \rangle $ 
is two-dimensional vector with integer components $k_1$, $k_2$
and, in particular, ${\bf 1} = \langle 1, 1 \rangle$.
We write $(m_1, \dots, m_k)$ for the greatest common factor of $m_1, \dots, m_k$.
As usual $\tau(k)$ is the number of positive divisors of $k$;
$r(k)$ is the number of solutions of the equation $m_1^2 + m_2^2 = k$ in integers $m_j$;
$\varphi(k)$ is the Euler function;
$\Omega(k)$ is the number of the prime factors of $k$, counted with the multiplicity;
$\chi(k)$ is the non-principal character modulo 4 and $L(s, \chi)$ is the corresponding Dirichlet's $L$-function.
We mark  by $\square$ an end of a proof, or its absence.

\section{Introduction and statement of the result.}

In 1937 Vinogradov \cite{Vin} considered the sum
\[
I^{(3)}(N) = \sum_{p_1 + p_2 + p_3 = N} (\log p_1)  ( \log p_2 )  ( \log p_3 ) .
\]
and proved that
\begin{equation} \label{10}
I^{(3)}(N) = \frac{1}{2} N^2 \mathfrak S^{(3)} (N) + O \left( N^2 \mathcal L^{-A} \right) ,
\end{equation}
where $A>0$ is an arbitrarily large constant and
\begin{equation} \label{11}
\mathfrak S^{(3)} (N) = \prod_{p \nmid N} \left( 1 + \frac{1}{(p-1)^3} \right) 
    \prod_{p \mid N} \left( 1 - \frac{1}{(p-1)^2} \right) .
\end{equation}

It is expected that a similar formula holds true for the sum
\[
  I^{(2)}(N) = \sum_{p_1 + p_2  = N} (\log p_1 ) ( \log p_2 ) ,
\]
but this has not been proved so far. However, using Vinogradov's method, one may establish that
$I^{(2)}(n)$ is close to $n \mathfrak S^{(2)} (n)$ for almost all $n \le N$. 
Here $\mathfrak S^{(2)} (n)$ is given by
\[
  \mathfrak S^{(2)} (n) 
  = \begin{cases}   c_0 \, \lambda(n) \quad & \text{for} \quad 2 \mid n , \\
                    0 & \text{for} \quad 2 \nmid n ,
         \end{cases}
\]
where
\begin{equation} \label{12}
  c_0  =   2 \prod_{p>2} \left( 1 - \frac{1}{(p-1)^2} \right)  , \qquad
\lambda(k) = \prod_{\substack{p \mid k \\ p > 2}} \frac{p-1}{p-2} .
\end{equation}
More precisely (see, for example, Vaughan \cite{Vaug}, Ch.2), for any constant $A>0$ we have
\begin{equation} \label{14}
   \sum_{n \le N } \left| I^{(2)}(n) -  n \mathfrak S^{(2)} (n) \right| 
   \ll N^2 \mathcal L^{-A} .
\end{equation}

Another classical achievement in prime number theory is the solution of the Hardy--Littlewood problem,
concerning the representation of large integers as a sum of two squares and a prime.
It was solved by Linnik~\cite{Linnik} and related problems have been studied
by Linnik, Hooley and other mathematicians
(see Hooley \cite{Hooley},~Ch.5). In particular, one can show that
\begin{equation} \label{20}
   \sum_{p \le N} r(p - 1) = 
   \pi N \mathcal L^{-1}
    \prod_{p > 2} \left( 1 + \frac{\chi(p)}{p(p-1)} \right)
    +    O \left( N \mathcal L^{-1-\theta_0} \left( \log \mathcal L \right)^5 \right) ,
\end{equation}
where
\begin{equation} \label{22}
\theta_0 = \frac{1}{2} - \frac{1}{4} e \log 2 = 0.0029\dots .
\end{equation}
A sharper estimate for the remainder term in \eqref{20} was established by Bredihin~\cite{Bredihin}.

In the present paper we prove a theorem which, in some sense, is a combination of
\eqref{10} and \eqref{20}.
Define
\begin{equation} \label{30}
R(N) = \sum_{p_1 + p_2 + p_3 = N} r(p_1 - 1) \; r(p_2 - 1) \; (\log p_1 ) (\log p_2 ) (\log p_3 ) 
\end{equation}
and
\begin{align} 
  \mathfrak S_{R} (N) 
      & =  
     \pi^2 \mathfrak S^{(3)}(N) 
      \prod_{p \nmid N(N-1)(N-2)} 
      \left( 1 + \chi(p) 
       \,\,
      \frac{2p^2 + p \chi(p) - 6p + 3 \chi(p)}{p^2(p^2 - 3 p + 3)} \right) 
      \notag \\
      & \qquad\qquad\qquad\;\;
      \times
      \prod_{p \mid N} 
       \qquad\;\,
      \left( 1 + \chi(p) \frac{2p^2 + p \chi(p) - 4p + 2\chi(p)}{p^2(p-1)(p-2)} \right) 
      \notag \\ 
      & \qquad\qquad\qquad\;
      \times
       \prod_{p \mid N-1} 
       \qquad
      \left( 1 + \chi(p) \frac{4p^2 - p \chi(p) - 6p + 3 \chi(p)}{p^2(p^2 - 3 p + 3)} \right) 
      \notag \\
      & \qquad\qquad\qquad\;
       \times
       \prod_{p \mid N-2} 
       \qquad
      \left( 1 + \chi(p) \frac{2p^2 - p^2 \chi(p) + p \chi(p) - 6p + 3 \chi(p)}{p^2(p^2 - 3 p + 3)} \right) ,
   \label{50}
\end{align}
where $\mathfrak S^{(3)}(N)$ is given by \eqref{11}.
\begin{theorem} \label{T1}
We have the following asymptotic formula
\begin{equation} \label{40}
  R(N) = \frac{1}{2} N^2 \mathfrak S_R(N) +  O \left(  N^2 \mathcal L^{- \theta_0}  \left( \log \mathcal L \right)^7 \right) ,
\end{equation}
where $\theta_0$ is the constant defined by \eqref{22}.
\end{theorem}

It is clear that if $2 \nmid N$ then $1 \ll \mathfrak S_{R} (N)  \ll 1$, so the main term in \eqref{40} dominates the remainder term provided that $N$ is a sufficiently large odd integer.

Theorem~\ref{T1} is related to a recent result of the author, which may be considered as a combination of \eqref{14}
and \eqref{20}. In the paper \cite{Tolev1} the sum
\[
 \sum_{p_1 + p_2 = n} r(p_1-1) \; (\log p_1 ) ( \log p_2 )
\]
was studied and it was proved that the expected asymptotic formula for it holds true
for almost all even integers $n \le N$. A similar problem was earlier considered by
Matom\"aki~\cite{Mato}.

The method used for the proof of Theorem~\ref{T1} can also be applied for 
finding asymptotic formulas for the sums
\[
\sum_{p_1 + p_2 + p_3 = N} r(p_1 - 1) \; \tau(p_2 - 1) \; (\log p_1) (\log p_2) (\log p_3 )
\]
and
\[
\sum_{p_1 + p_2 + p_3 = N} \tau(p_1 - 1) \; \tau(p_2 - 1) \; (\log p_1) (\log p_2) (\log p_3) .
\]

It would be interesting to consider the ternary Goldbach equation with weights 
of the above type attached to all of the variables. 
We would be in a position to attack this problem if we had more information about the number of solutions of the ternary equation with all prime variables lying in independent arithmetic progressions with large moduli. 
However, the best result of this type available in the literature at present, which is due to the author
\cite{Tolev2} and improves a theorem of K.Halupczok~\cite{Halup}, is not strong enough for our aims.

\section{Some lemmas.}

First we consider the Goldbach binary problem with one prime variable
lying in a given interval and belonging to an arithmetic progression. Suppose that $n \le N$, let $k$ and $l$ be integers with
$(k,l)=1$ and let $J \in \mathcal J$. We denote
\begin{align} 
I^{(2)}_{k, l}(n , J) 
  & = 
  \sum_{\substack{p_1 + p_2 = n \\ p_1 \equiv l \pmod{k} \\ p_1 \in J}} 
  (\log p_1) ( \log p_2 ) ; 
  \label{52} \\
  & \notag \\
\mathfrak S^{(2)}_{k, l}(n)
  & = 
 \begin{cases}
   c_0 
   \lambda(nk)
       \qquad
      &    \text{if} \;\; (k, n-l) = 1 \; \; \text{and} \;\;2 \mid n , \\
       0 
       &  \text{otherwise} ;
      \end{cases} 
          \label{54} \\
          & \notag \\
  \Phi^{(2)}(n , J) 
   & = 
   \sum_{\substack{m_1 + m_2 = n \\ m_1 \in J}} 1 ;
   \label{55} \\
          & \notag \\
  \Delta^{(2)}_{k, l} (n, J) 
    & = 
    I^{(2)}_{k, l}(n , J) -
   \frac{\mathfrak S^{(2)}_{k, l}(n)}{\varphi(k)} \Phi^{(2)}(n , J) .
  \label{56}
\end{align}
If $J = [1, N]$ then we write for simplicity $I^{(2)}_{k, l}(n)$, 
$ \Phi^{(2)}(n) $ $ (= n-1 $) and $\Delta^{(2)}_{k, l} (n)$.

Our first lemma is a generalization of \eqref{14} and states that  $  \Delta^{(2)}_{k, l} (n, J) $ 
is small on average with respect to $k$ and $n$ and uniformly for $l$ and $J$.
More precisely, we have
\begin{lemma} \label{L1}
For any constant $A>0$ there exist $B=B(A) > 0$ such that
\[
\sum_{k \le \sqrt{N} \mathcal L^{-B} } \max_{(l, k)=1} \max_{J \in \mathcal J}
  \sum_{ n \le N } \left|
   \Delta^{(2)}_{k, l} (n, J)
     \right|
     \ll N^2 \mathcal L^{-A} .
\]
\end{lemma}

This lemma is very similar to results of Mikawa~\cite{Mikawa} and Laporta~\cite{Laporta}. 
These authors study the equation $p_1 - p_2 = n$ and without the condition $p_1 \in J$. 
However, inspecting the arguments presented in  \cite{Laporta}, 
the reader will readily see that the proof of Lemma~\ref{L1} can be obtained is the same manner. $\square$

Next we consider Goldbach's ternary problem with two primes from arithmetic progressions and belonging to given 
intervals. Suppose that ${\bf k} = \langle k_1, k_2 \rangle$ and ${\bf l} = \langle l_1, l_2 \rangle$
are two-dimensional vectors with integer components and let
${\bf J} = \langle J_1, J_2 \rangle$ be a pair of intervals $J_1, J_2 \in \mathcal J$.
We denote
\begin{align} 
  I^{(3)}_{ {\bf k} , \, {\bf l}}(N, {\bf J}) 
   &  = 
  \sum_{\substack{p_1 + p_2 + p_3 = N \\ p_i \equiv l_i \pmod{k_i} \\ p_i \in J_i \, , \; i=1,2}} 
  (\log p_1 ) (\log p_2 ) (\log p_3 ) , 
    \label{60} \\
    & \notag \\
  \Phi^{(3)} (N, {\bf J} ) 
   & = 
   \sum_{\substack{m_1 + m_2 + m_3 = N \\ m_i \in J_i \, , \; i=1,2 }} 1 .
   \notag
\end{align}

Using the notations of K.Halupczok~\cite{Halup},
we define $\mathfrak S^{(3)}_{{\bf k} , {\bf l}}(N) $ in the following way.
Consider the sets of primes
\begin{align}
  \mathcal A 
    & = \{ p  \; : \; p \nmid k_1 k_2 , \; p \mid N \} ; 
      \notag \\
   \mathcal B 
     & = 
     \{ p \; : \; p \nmid k_1 k_2 N \} ;
     \notag \\
   \mathcal C 
    & = \{ p \; : \; p \mid k_1 , \; p \nmid  k_2 , \; p \mid N - l_1 \} \cup \{ p \; : \;
               p \mid k_2 , \; p \nmid  k_1 , \; p \mid N - l_2 \} ;
            \notag \\
   \mathcal D 
    & = \{ p \; : \; p \mid k_1 , \; p \nmid  k_2 , \; p \nmid N - l_1 \} \cup \{ p \; : \;
               p \mid k_2 , \; p \nmid  k_1 , \; p \nmid N - l_2 \};
            \notag \\
   \mathcal E
    & = 
    \{ p \; : \; p \mid k_1 , \; p \mid  k_2 , \; p \mid N - l_1  - l_2 \};
   \notag \\
   \mathcal F
    & = 
    \{ p \; : \; p \mid k_1 , \; p \mid  k_2 , \; p \nmid N - l_1  - l_2 \}.
   \notag
\end{align}        
If $\mathcal E \not= \emptyset$ then we assume that
\begin{equation}   \label{69} 
   \mathfrak S^{(3)}_{{\bf k} , {\bf l}}(N)  = 0 .
\end{equation}
If $\mathcal E = \emptyset$ then we put
\begin{equation}   \label{70} 
  \mathfrak S^{(3)}_{{\bf k} , {\bf l}}(N) 
    = 
      \prod_{p \in \mathcal A \cup \mathcal D } \left( 1 - \frac{1}{(p-1)^2} \right) 
   \prod_{p \in \mathcal B } \left( 1 + \frac{1}{(p-1)^3} \right) 
   \prod_{p \in  \mathcal C \cup \mathcal F} \left( 1 + \frac{1}{p-1} \right) .
\end{equation}

We also define
\begin{equation}   \label{75} 
  \Delta^{(3)} _{ {\bf k} , \, {\bf l}}(N , {\bf J})  
    =
        I^{(3)}_{ {\bf k} , \, {\bf l}}(N , {\bf J})  -
  \frac{\mathfrak S^{(3)}_{{\bf k} , {\bf l}}(N) \; \Phi^{(3)} (N ,  {\bf J}) }{\varphi(k_1) \varphi(k_2)} .
\end{equation}
If $J_1 = J_2 = [1, N]$ then we write for simplicity
$   I^{(3)}_{{\bf k} , \, {\bf l}}(N)  $, $\Phi^{(3)}(N)$ ($ = N^2/2 + O (N)$)
and $\Delta^{(3)} _{ {\bf k} , \, {\bf l}}(N)  $.

The next lemma is analogous to Lemma~\ref{L1}
and states that $\Delta^{(3)} _{ {\bf k} ,  {\bf l}}(N , {\bf J})  $
is small on average with respect to ${\bf k}$ and uniformly for ${\bf l}$ and ${\bf J}$. 
More precisely, we have
\begin{lemma} \label{L2}
For any constant $A>0$ there exist $C=C(A) > 0$ such that
\[
\sum_{k_1 \le \sqrt{N} \mathcal L^{-C} }
\sum_{ k_2 \le \sqrt{N} \mathcal L^{-C} } \; \max_{\substack{ (l_i, k_i)=1 \\ i=1,2}} 
  \max_{\substack{ J_i \in \mathcal J \\ i=1,2 }}
    \; \left|  
    \Delta^{(3)} _{ {\bf k} ,  {\bf l}}(N , {\bf J})  
     \right| \;
     \ll N^2 \mathcal L^{-A} .
\]
\end{lemma}

This statement is 
slightly more general than a theorem 
from author's recent paper \cite{Tolev2}, which improves a result of K.Halupczok \cite{Halup}.
There are no conditions $p_i \in J_i$ in the theorems of \cite{Tolev2} and \cite{Halup},
but the reader can easily verify that the methods developed in these articles imply also the validity of Lemma~\ref{L2}. $\square$

In several occasions we will need the following simple
\begin{lemma}  \label{L3}                  
Suppose that $j \in \{1, -1 \}$ and let $m$, $k$, $l$, $n$ be natural numbers.
Then the quantities
$\mathfrak S^{(2)}_{ 4m ,  1 + jm }(n) $ 
and
$\mathfrak S^{(3)}_{\langle k, 4m \rangle , \langle l, 1 + jm \rangle } (n) $
do not depend on $j$.
\end{lemma}  
The proof follows directly from the definitions of $\mathfrak S^{(2)}_{k, l}(n)$ and $\mathfrak S^{(3)}_{{\bf k} , {\bf l}}(n)$.
We leave the easy verification to the reader.
$\square$

The next two lemmas are due to C.Hooley and play an essential role in the proof of \eqref{20}, 
as well as in the solutions of other related problems. 
\begin{lemma} \label{L4}
For any constant $\omega > 0 $ we have
\[
\sum_{p \le N} \left| 
   \sum_{\substack{ d \mid p - 1 \\\sqrt{N} \mathcal L^{-\omega} < d <  \sqrt{N} \mathcal L^{\omega} } }
   \chi(d)
\right|
   \ll N \mathcal L^{-1 - \theta_0} \left( \log \mathcal L \right)^{5} ,
\]
where $\theta_0$ is defined by \eqref{22}. The constant in the Vinogradov symbol depends only on $\omega$.
\end{lemma}

\begin{lemma} \label{L5}
For any constant $\omega > 0 $ we have
\[
\sum_{p \le N} \left| 
   \sum_{\substack{ d \mid p - 1 \\\sqrt{N} \mathcal L^{-\omega} < d <  \sqrt{N} \mathcal L^{\omega} } }
   \chi(d)
\right|^2
   \ll N \mathcal L^{-1} \left( \log \mathcal L \right)^7 ,
\]
where the constant in the Vinogradov symbol depends on $\omega$.
\end{lemma}

The proofs of very similar results (with $\omega = 48$ and with the condition $d \mid N - p$ 
rather than $d \mid p-1$) are available in \cite{Hooley},~Ch.5. The reader will easily see that 
the methods used there yield also the validity of Lemma~\ref{L4} and Lemma~\ref{L5}.
$\square$

The next lemma is analogous to another result of Hooley from \cite{Hooley},~Ch.5.
\begin{lemma} \label{L6}
Let $n$ be an integer satisfying $1 \le n \le N$. Suppose 
that $\omega > 0 $ is a constant and let $\mathcal P = \mathcal P_{\omega}(N)$ 
be the set of primes $p \le N$ such that
$p-1$ has a divisor lying between $ \sqrt{N} \mathcal L^{-\omega} $ and $ \sqrt{N} \mathcal L^{\omega} $.
Then we have
\begin{equation} \label{85}
 \sum_{\substack{p_1 + p_2 = n \\ p_1 \in \mathcal P}} 1
 \ll N \mathcal L^{-2 - 2 \theta_0} \left( \log \mathcal L \right)^6 ,
\end{equation}
where $\theta_0$ is defined by \eqref{22}.
The constant in the Vinogradov symbol depends only on $\omega$.
\end{lemma}

{\bf Proof:} We proceed as in \cite{Hooley}, Ch.~5, Sec.~7. Denote the sum in the left side of 
\eqref{85} by $\Sigma$ and let $D_1 = \sqrt{N} \mathcal L^{- \omega}$, $D_2 = \sqrt{N} \mathcal L^{ \omega}$. Suppose that $\alpha$ is a real number satisfying
$1 < \alpha < 3/2$. We have
\begin{equation} \label{87}
  \Sigma \le 
   \sum_{\substack{p_1 + p_2 = n \\ \Omega(p_1 - 1) \le \alpha \log \mathcal L \\p_1 \in \mathcal P}} 1 + 
   \sum_{\substack{p_1 + p_2 = n \\ \Omega(p_1 - 1) > \alpha \log \mathcal L }} 1 =
  \Sigma_1 + \Sigma_2 ,
\end{equation}
say. 

Consider first $\Sigma_1$. We have
\[
 \Sigma_1 \le  \sum_{\substack{p_1 + p_2 = n \\ \Omega(p_1 - 1) \le \alpha \log \mathcal L }} 
  \sum_{\substack{d \mid p_1 - 1 \\ D_1 < d < D_2 }} 1 =
  \sum_{p_1 + p_2 = n} \sum_{\substack{md = p_1 - 1 \\ D_1 < d < D_2 \\ \Omega(md) \le \alpha \log \mathcal L }} 1 .
\]
Obviously, the inequality  $\Omega(md) \le \beta$ implies the validity of at least one of the 
inequalities
$\Omega(m) \le \frac{1}{2} \beta$ and  $\Omega(d) \le \frac{1}{2} \beta$. Hence
\begin{equation} \label{89}
 \Sigma_1 \le 
  \sum_{p_1 + p_2 = n} \sum_{\substack{md = p_1 - 1 \\ D_1 < d < D_2 \\ \Omega(d) \le \frac{1}{2} \alpha \log \mathcal L }} 1 
  +  \sum_{p_1 + p_2 = n} \sum_{\substack{md = p_1 - 1 \\ D_1 < d < D_2 \\ \Omega(m) \le \frac{1}{2} \alpha \log \mathcal L }} 1 
   = \Sigma_1^{(1)} +  \Sigma_1^{(2)} ,
\end{equation}
say. Consider $\Sigma_1^{(2)}$. 
The conditions imposed in its definition imply $m < D_2$ and clearly
\begin{equation} \label{90}
 \Sigma_1^{(2)} \ll 
  \sum_{p_1 + p_2 = n} \sum_{\substack{md = p_1 - 1 \\ D_1 \mathcal L^{-10} < m < D_2 \\ \Omega(m) \le \frac{1}{2} \alpha \log     \mathcal L }} 1 +  N \mathcal L^{-10 }
   =
 \Sigma_1^{(3)} +  N \mathcal L^{-10 } ,
\end{equation}
say. Obviously $\Sigma_1^{(1)} \le \Sigma_1^{(3)}$ and from this inequality,
\eqref{89} and \eqref{90} we find
\begin{equation} \label{91}
\Sigma_1  \ll \Sigma_1^{(3)} + N \mathcal L^{-10} .
\end{equation}

To estimate $ \Sigma_1^{(3)}$ we change the order of summation and find
\begin{equation} \label{92}
 \Sigma_1^{(3)}  = 
 \sum_{\substack{ D_1 \mathcal L^{-10} < m < D_2 \\ \Omega(m) \le \frac{1}{2} \alpha \log \mathcal L}}
 \lambda_m(n) , 
\end{equation}
where $\lambda_m(n)$ is the number of primes $p < n$
satisfying $p \equiv 1 \pmod{m}$ and such that
$n-p$ is a prime too.
We apply Theorem~2.4 of Halberstam and Richert \cite{H-R} (with $x=y=n$, $k=m$, $l=1$, $g=1$, $a_1=-1$, $b_1=n$) and find
\[
 \lambda_m(n) \ll \prod_{p \mid mn} \left( 1 - \frac{1}{p} \right)^{-1} \cdot \frac{m}{\varphi(m)} \cdot
 \frac{n/m}{ \log^2 (n/m)  } \ll \frac{N}{m} \, \mathcal L^{-2} \, (\log \mathcal L)^2 .
\]
We substitute this upper bound for $\lambda_m(n)$ in \eqref{92} and then proceed precisely  as in \cite{Hooley}, Ch.5, Sec.~7
to find
\begin{equation} \label{94}
 \Sigma_1^{(3)} \ll N \mathcal L^{-2} (\log \mathcal L)^2 \sum_{\substack{D_1 \mathcal L^{-10}< m < D_2 \\ 
 \Omega(m) \le \frac{ \alpha }{ 2 } \log \mathcal L}} \frac{1}{m} \ll
 N \mathcal L^{\gamma(\alpha/2) - 3} (\log \mathcal L)^3 ,
\end{equation}
where 
\begin{equation} \label{94.5}
    \gamma(c) = c - c \log c .
\end{equation}
From \eqref{91} and \eqref{94} we get
\begin{equation} \label{95}
 \Sigma_1 \ll  N \mathcal L^{\gamma(\alpha/2) - 3} (\log \mathcal L)^3 .
\end{equation}

Consider now $\Sigma_2$.  We have
\begin{equation} \label{96}
 \Sigma_2 \le 
    \sum_{\substack{p_1 + p_2 = n \\ \alpha \log \mathcal L < \Omega(p_1 - 1) \le 10 \log \mathcal L }} 1 +
    \sum_{\substack{m \le N \\ \Omega(m) > 10 \log \mathcal L}} 1 = 
    \Sigma_2^{(1)} + \Sigma_2^{(2)} ,
\end{equation}
say. It is shown in \cite{Hooley}, Ch.5, Sec.~7 that
\begin{equation} \label{98}
\Sigma_2^{(2)} \ll N \mathcal L^{-4} .
\end{equation}

Consider $\Sigma_2^{(1)}$.
Denote by $\mathcal R_N$ the set of integers $m \le N$ composed only of primes 
$\le N^{\frac{1}{20 \log \mathcal L }}$.
Applying the method of \cite{Hooley}, Ch.5, Sec.~7
we get
\begin{align}
 \Sigma_2^{(1)} 
   & \le
   \sum_{\substack{p_1 + p_2 = n \\ \alpha \log \mathcal L < \Omega(p_1 - 1)  \\
   p_1 - 1 \not\in \mathcal R_N}} 1 + 
   \sum_{\substack{p_1 + p_2 = n \\  \Omega(p_1 - 1) \le 10 \log \mathcal L \\
   p_1 - 1 \in \mathcal R_N}} 1 
   \notag \\
   & \notag \\
   & = 
    \sum_{\substack{p_1 + p_2 = n \\ \alpha \log \mathcal L < \Omega(p_1 - 1)  \\
   p_1 - 1 \not\in \mathcal R_N}} 1 
      +  O (\sqrt{N}) 
      \notag \\
      & \notag \\
    & =
     \Sigma_2^{(3)} +   O (\sqrt{N}) ,
   \label{100}
\end{align}
say.  If $p_1 - 1 \not\in \mathcal R_N$ then there exists a prime $q > N^{\frac{1}{20 \log \mathcal L}}$ such that
$p_1 - 1 = q r$ for some positive integer $r < N^{1 - \frac{1}{20 \log \mathcal L}}$. On the other hand, from the condition
$ \Omega(p_1 - 1) >  \alpha \log \mathcal L  $ it follows that $ \Omega(r) > \alpha \log \mathcal L - 1 $.
Therefore
\begin{equation} \label{102}
   \Sigma_2^{(3)} \le \sum_{\substack{ r <  N^{1 - \frac{1}{20 \log \mathcal L}} \\
   \Omega(r) >  \alpha \log \mathcal L - 1     }} \varkappa_r(n) ,
\end{equation}
where $\varkappa_r(n)$ is the number of primes $q \le N/r$ such that $r q + 1 $ and $-r q + n -1 $ are primes too.
We apply again Theorem~2.4 of \cite{H-R} (with $x = y = N/r$, $k=1$, $g=2$, $a_1 = r$, $a_2= -r$, $b_1=1$, $b_2 = n-1$,
$E= - r^3 n$) to find
\begin{equation} \label{103}
   \varkappa_r(n) \ll \prod_{p \mid E} \left( 1 - \frac{1}{p} \right)^{\rho(p) - 2 } \cdot
                     \prod_{p \mid n-1} \left( 1 - \frac{1}{p} \right)^{-1} \cdot
                     \frac{N/r}{\log^3 (N/r)} ,
\end{equation}
where $\rho(p)$ is the number of solutions of the congruence
$(a_1 m + b_1)(a_2 m + b_2) \equiv 0 \pmod{p}$. It is easy to verify that
\begin{equation} \label{106}
\rho(p) = \begin{cases}
    p \quad & \text{for} \quad p \mid r, \; p \mid n-1  ,\\
    0 \quad & \text{for} \quad p \mid r, \; p \nmid n-1  ,\\
    1 \quad & \text{for} \quad p \nmid r, \; p \mid n-2 , \\
    2 \quad & \text{for} \quad p \nmid r, \; p \nmid n-2  .
     \end{cases}
\end{equation}
From \eqref{103} and \eqref{106} it follows that
\[
  \varkappa_r(n) \ll \frac{N}{r \log^3 (N/r)} (\log \mathcal L)^3 .
\]
We substitute this upper bound for $\varkappa_r(n)$ in \eqref{102} and then
we notice that the inequality
$r < N^{1 - \frac{1}{20 \log \mathcal L}}$ 
implies
$\log (N/r) \gg \mathcal L (\log \mathcal L)^{-1}$. 
We find
\[
  \Sigma_2^{(3)} \ll
     N \mathcal L^{-3} (\log \mathcal L)^6
   \sum_{\substack{ r <  N \\
   \Omega(r) > \alpha \log \mathcal L - 1     }} \frac{1}{r} .
\]
Now we estimate the sum over $r$ in the way proposed in \cite{Hooley}, Ch.~5, Sec.~7 to get
\begin{equation} \label{110}
  \Sigma_2^{(3)} \ll
     N \mathcal L^{\gamma(\alpha) - 3} (\log \mathcal L)^6 ,
\end{equation}
where $\gamma(c)$ is defined by \eqref{94.5}.
Using \eqref{96}, \eqref{98}, \eqref{100}  and \eqref{110}
we obtain
\begin{equation} \label{104}
 \Sigma_2 
\ll
     N \mathcal L^{\gamma(\alpha) - 3} (\log \mathcal L)^6 .
\end{equation}

From \eqref{87}, \eqref{95} and \eqref{104} it follows that
\[
 \Sigma \ll N \left( \mathcal L^{\gamma(\alpha/2)- 3}  + \mathcal L^{\gamma(\alpha) - 3} \right) (\log \mathcal L)^6 .
\]
We choose $\alpha$ from the condition
$\gamma(\alpha / 2 ) = \gamma(\alpha)$, which gives
$\alpha = e /2$. This completes the proof of the lemma. $\square$

\section{Proof of Theorem~\ref{T1}.}

\paragraph{Beginning.}

We put
\begin{equation} \label{150}
  D = \sqrt{N} \, \mathcal L^{ - B(10) - C(10) - 1   } ,
\end{equation}
where $B(A)$ and $C(A)$ are specified respectively in Lemma~\ref{L1} and Lemma~\ref{L2}.
Obviously
\begin{equation} \label{160}
   r(m)  = 4 \sum_{d \mid m} \chi(d) = 4 \left( r_1 (m) + r_2 (m) + r_3 (m)  \right) ,
\end{equation}
where
\begin{equation} \label{171}
r_1 (m) =  \sum_{\substack{ d \mid m \\ d \le D} } \chi(d) , \qquad
r_2 (m) =  \sum_{\substack{ d \mid m \\ D < d < N/D} } \chi(d) , \qquad
r_3 (m) =  \sum_{\substack{ d \mid m \\ d \ge N/D} } \chi(d) .
\end{equation} 
Hence using \eqref{30} and \eqref{160} we get
\begin{equation} \label{172}
 R(N) = 16  \sum_{1 \le i, j \le 3} \;
   R_{i, j}(N) \; ,
\end{equation}  
where
\begin{equation} \label{173}
 R_{i, j} (N) = 
  \sum_{p_1 + p_2 + p_3 = N} 
    r_i(p_1 - 1) \; r_j(p_2 - 1) \;  (\log p_1) ( \log p_2 ) ( \log p_3 ) .
\end{equation}  
We shall prove that the main term 
in \eqref{40} comes from $R_{1,1}(N)$ and the other sums $R_{i, j}(N)$ contribute only to the remainder
term. Because of the symmetry we have to consider only the expressions $R_{i, j}(N)$ with $i \le j$.
  
\paragraph{The evaluation of $R_{1, 1}(N)$.} 

Using \eqref{60}, \eqref{75}, 
\eqref{171} and \eqref{173} we get
\begin{equation} \label{173.1}
  R_{1, 1}(N)  
  =     \sum_{d_1, d_2 \le D} \chi(d_1) \chi(d_2) \; 
    I^{(3)}_{{\bf d}, {\bf 1}} (N)
    = 
    R_{1, 1}'(N) +  R_{1, 1}^*(N)  ,
\end{equation}  
where
\begin{align}
    R_{1, 1}'(N)
    & =
    \sum_{d_1, d_2 \le D} \frac{\chi(d_1) \chi(d_2)}{\varphi(d_1) \varphi(d_2)}
    \mathfrak S_{{\bf d}, {\bf 1}}^{(3)}(N) \Phi^{(3)}(N) ,
     \notag \\
     & \notag \\
     R_{1, 1}^*(N)
    & =
    \sum_{d_1, d_2 \le D} \chi(d_1) \chi(d_2) \; 
    \Delta_{{\bf d}, {\bf 1}}^{(3)}(N) .
     \notag
\end{align}     
From \eqref{150} and Lemma~\ref{L2} it follows that
\begin{equation} \label{173.7}
     R_{1, 1}^*(N) \ll N^2 \mathcal L^{-1} .
\end{equation}  
Consider $R_{1, 1}'(N)$. It is clear that
\begin{equation} \label{174}
     R_{1, 1}'(N)  =  \frac{1}{2}N^2 \; \Gamma(N) + O \left( N^{1 + \varepsilon}  \right) ,
\end{equation}  
where
\begin{equation} \label{175}
  \Gamma (N) =  \sum_{d_1, d_2 \le D} \frac{\chi(d_1) \chi(d_2)}{\varphi(d_1) \varphi(d_2)}
    \mathfrak S_{{\bf d}, {\bf 1}}^{(3)}(N) . 
\end{equation}  

It remains to establish an asymptotic formula for $\Gamma (N)$. 
The calculations are long and complicated, but rather routine and straightforward. 
We point out only the main steps and leave the details to the reader.

Using \eqref{69}, \eqref{70} and \eqref{175} we find
\begin{equation} \label{176}
  \Gamma (N) = \sum_{d \le D} \frac{\chi(d)}{\varphi(d)} \psi_N(d)
    \sum_{\substack{t \le D \\ (t, d, N-2) = 1}} f_{N, d}(t) ,  
\end{equation}  
where
\begin{align}
  \psi_{N}(d) 
    & = 
    \prod_{p \nmid dN} \left( 1 + \frac{1}{(p-1)^3} \right)
    \prod_{ p \in \mathcal U_1 \cup \; \mathcal U_2} \left( 1 - \frac{1}{(p-1)^2} \right)
    \prod_{ \substack{ p \mid N-1 \\ p \mid d }} \left( 1 + \frac{1}{p-1} \right) ,
    \label{177} \\
    & \notag \\
   f_{N, d} (t)
    & =
     \frac{\chi(t)}{\varphi(t)}
     \prod_{ \substack{ p \nmid d N \\ p \mid t }} \left( 1 + \frac{1}{(p-1)^3} \right)^{-1} 
     \prod_{\substack{ p \in \mathcal U_1 \cup \; \mathcal U_2 \\ p \mid t}} \left( 1 - \frac{1}{(p-1)^2} \right)^{-1}
     \prod_{ p \mid (d, t,  N-1)} \left( 1 + \frac{1}{p-1} \right)^{-1} 
     \notag \\
      & \notag \\
     & \qquad\qquad
      \times
      \prod_{ \substack{ p \nmid d (N-1) \\ p \mid t }} \left( 1 - \frac{1}{(p-1)^2} \right)
      \prod_{\substack{ p \in \mathcal U_3 \cup \; \mathcal U_4 \\ p \mid t}} \left( 1 + \frac{1}{p-1} \right)
    \label{178} 
\end{align}
and 
\begin{align}
  & \mathcal U_1 = 
  \{ p \; : \; p \mid N , \; p \nmid d \} , 
   & \mathcal 
    U_2   = 
  \{ p \; : \; p \nmid N-1 , \; p \mid d \} ,
    \notag  \\
   & \mathcal U_3 = 
  \{ p \; : \; p \mid N-1 , \; p \nmid d \} , 
    & \mathcal  U_4 = 
  \{ p \; : \; p \nmid N-2 , \; p \mid d \} .
    \notag
\end{align}

First we evaluate the sum over $t$ in \eqref{176}. From \eqref{178} it follows that 
\begin{equation} \label{179}
  f_{N, d}(t) \ll (\log \mathcal L)^2  \;  t^{-1} \, \log \log (10 t)
\end{equation}
with absolute constant in the Vinogradov symbol.
Hence the corresponding Dirichlet series
\[
  F_{N, d} (s) = \sum_{\substack{t = 1 \\ (t, d, N-2)=1}}^{\infty} f_{N, d}(t) \, t^{-s}
\]
is absolutely convergent in $Re \; (s) > 0$. Clearly $f_{N, d}(t)$ is multiplicative with respect to $t$ and
applying Euler's identity we find
\[
 F_{N, d} (s) = \prod_{p \nmid (d, N-2)} T_{N, d}(p, s) , \qquad
 T_{N, d}(p , s) = 1 + \sum_{l=1}^{\infty} f_{N, d} \left( p^l \right) p^{-ls} .
\]
From \eqref{178} we establish that
\[
 T_{N, d}(p , s) = \left( 1 - \chi(p) p^{-s-1} \right)^{-1} 
 \left( 1 + \chi(p) p^{-s-1} Y_{N, d} (p) \right) , 
\]
where
\begin{equation} \label{181}
  Y_{N, d} (p) =
  \begin{cases}
    (p-3) (p^2 - 3 p + 3)^{-1}  & \text{if} \qquad p \nmid dN(N-1)(N-2) , \\
    2(p-2)^{-1}  & \text{if} \qquad  p \mid d , \; p \nmid N(N-1)(N-2) , \\
    (p-1)^{-1}  & \text{if}  \qquad p \nmid d , \; p \mid N , \\
    2 (p-2)^{-1}  & \text{if}  \qquad p \mid d , \; p \mid N , \\
    (2p-3)(p^2 - 3 p + 3)^{-1}  & \text{if}  \qquad p \nmid d , \; p \mid N-1 , \\
    (p-1)^{-1}  & \text{if}  \qquad p \mid d , \; p \mid N-1 , \\
    (p-3)(p^2 - 3 p + 3)^{-1}  & \text{if}  \qquad p \nmid d , \; p \mid N-2 . \\
  \end{cases}
\end{equation}

Hence we get
\begin{equation} \label{182}
 F_{N, d}(s) = L(s+1, \chi) \prod_{p \mid (d, N-2)}  \left( 1 - \chi(p) p^{-s-1} \right) 
    \prod_{p \nmid (d, N-2)}
   \left( 1 + \chi(p) p^{-s-1} Y_{N, d} (p) \right) .
\end{equation}
From this formula it follows that $F_{N, d}(s)$ has an analytic continuation to $Re \; (s) > -1$.
Using \eqref{181} and the simplest bound for $L(s+1, \chi)$ we get
\begin{equation} \label{183}
 F_{N, d} (s) \ll N^{\varepsilon} \, T^{1/6} \qquad \text{for} \qquad Re \; (s) \ge -\frac{1}{2} , \quad
    |Im \; (s)| \le T .
\end{equation}

We apply the version of Perron's formula given at Tenenbaum~\cite{Tenen}, Ch. II.2 and also
\eqref{179} to find
\[
  \sum_{\substack{t \le D \\ (t, d, N-2)= 1 }} f_{N, d}(t) =
  \frac{1}{2 \pi i} \int_{\varkappa - i T}^{\varkappa + i T} F_{N, d} (s) \frac{D^s}{s} \, d s 
  +  O \left( \sum_{t=1}^{\infty} \frac{N^{\varepsilon} \, D^{\varkappa} \, \log \log (10 t) }
  {t^{1 + \varkappa} \left( 1 + T \left| \log \frac{D}{t} \right| \right) } \right) ,
\]
where $\varkappa = 1/100$, $T = N^{3/4} $. It is easy to verify that the remainder term above is
$ O \left(  N^{-1/100} \right) $ and applying the residue theorem we see that the main term is equal to
\[
  F_{N, d}(0) + \frac{1}{2 \pi i} \left(  \int_{\varkappa - i T}^{-1/2 - i T} +
    \int_{ -1/2 - i T}^{ -1/2 + i T} + \int_{ -1/2 + i T}^{\varkappa + i T}\right)  F_{N, d} (s) \frac{D^s}{s} \, d s .
\]
From \eqref{183} it follows that the contribution from the above integrals is $ O \left(  N^{-1/100} \right) $.
Hence 
\begin{equation} \label{185}
 \sum_{\substack{t \le D \\ (t, d, N-2)= 1 }} f_{N, d}(t) =   F_{N, d}(0) + O \left( N^{-1/100} \right) .
\end{equation}
Obviously, using \eqref{182} we get
\begin{equation} \label{186}
 F_{N, d}(0) = \frac{\pi}{4} \prod_{p \mid (d, N-2)}  \left( 1 - \frac{\chi(p)}{p} \right) 
    \prod_{p \nmid (d, N-2)}
   \left( 1 + \frac{\chi(p)}{p} Y_{N, d} (p) \right) .
\end{equation}

We use \eqref{176}, \eqref{177}, \eqref{181}, \eqref{185} and \eqref{186} to find a new expression for 
$\Gamma(N)$ and after some calculations we obtain
\begin{equation} \label{187}
  \Gamma(N) = \frac{\pi}{4} \, \mathfrak S^{(3)}(N) \, \Xi (N) \sum_{d \le D} g_N(d) 
     + O \left( N^{\varepsilon - 1/100} \right) ,
\end{equation}
where $\mathfrak S^{(3)}(N)$ is defined by \eqref{11},
\begin{align} 
 \Xi(N) 
   & = 
   \prod_{p \nmid N(N-1)} \left( 1 + \frac{\chi(p)(p-3)}{p(p^2- 3 p + 3)} \right) 
          \prod_{p \mid N} \left( 1 + \frac{\chi(p)}{p(p-1)} \right) 
        \notag \\
        & \notag \\
   & \qquad \qquad \qquad
      \times       
          \prod_{p \mid N-1} \left( 1 + \frac{\chi(p)(2p-3)}{p(p^2- 3 p + 3)} \right) ,
  \label{189}
\end{align}
and
\begin{align}
  g_N(d) 
    & = 
      \frac{\chi(d)}{\varphi(d)}
      \prod_{\substack{p \mid d \\   p \nmid N(N-1)(N-2)}}
      \frac{1 + \frac{2 \chi(p)}{p(p-2)} }  {1 + \frac{\chi(p) (p-3)}{p(p^2 - 3 p + 3)} } 
        \;\;
      \prod_{p \mid (d, N) }
      \frac{1 + \frac{2 \chi(p)}{p(p-2)} }  {\left( 1 - \frac{1}{(p-1)^2}  \right) \left( 1 + \frac{\chi(p)}{p(p-1)} \right) }
      \notag \\
      & \notag \\
      & \qquad \qquad
        \times
        \prod_{\substack{p \mid d \\   p \nmid N}} 
          \left( 1 + \frac{1}{(p-1)^3} \right)^{-1} 
          \;\;\,
          \prod_{p \mid (d, N-1) }
      \frac{\left( 1 + \frac{1}{p-1} \right) \left( 1 + \frac{\chi(p)}{p(p-1)} \right) }  
       { 1  + \frac{\chi(p)(2 p - 3 )}{p(p^2-3p + 3)}  }
          \notag \\
      & \notag \\
      & \qquad \qquad
        \times
          \prod_{\substack{p \mid d \\   p \nmid N - 1 }} 
          \left( 1 - \frac{1}{(p-1)^2} \right)
              \prod_{p \mid (d, N-2) }
      \frac{ 1 - \frac{\chi(p)}{p} }  
       { 1  + \frac{\chi(p)( p - 3 )}{p(p^2-3p + 3)}  } .
       \label{190}
\end{align}

It is clear that $g_N(d)$ is multiplicative with respect to $d$ and satisfies
\[
   g_N(d) \ll (\log \log (10 d) )^3 \; d^{-1} ,
\]
where the constant in Vinogradov's symbol is absolute.
Hence the Dirichlet series 
\[
  G_N(s) = \sum_{d=1}^{\infty} g_N(d) \, d^{-s}
\]
is absolutely convergent in $Re \; (s) > 0$ and applying the Euler identity we get
\begin{equation} \label{192}
  G_N(s) = \prod_p H_N(p, s) , \qquad H_N(p, s) = 1 + \sum_{l=1}^{\infty} g_N \left( p^l \right) p^{-ls}.
\end{equation}
From \eqref{190} and \eqref{192} we find
\[
   H_N(p, s) = \left( 1 - \chi(p) p^{-s-1} \right)^{-1} 
   \left( 1 + \chi(p) p^{-s-1} K_N(p) \right) ,
\]
where
\begin{equation} \label{196}
  K_N(p) = 
    \begin{cases}
     \frac{p^2 + p \chi(p) - 3 p + 3 \chi(p)}{p^3- 3 p^2 + 3 p + p \chi(p) - 3 \chi(p)} 
         \quad & \text{if} \quad p \nmid N(N-1)(N-2) , \\
          & \\
         \frac{p^2 + p \chi(p)- 2p + 2 \chi(p)}{p^3-3p^2 + p \chi(p) + 2p - 2 \chi(p)} 
         \quad & \text{if} \quad p \mid N , \\
         & \\
       \frac{2p^2 - 3 p - p \chi(p) + 3 \chi(p)}{p^3- 3 p^2 + 3 p + 2 p \chi(p) - 3 \chi(p)} 
         \quad & \text{if} \quad p \mid N-1 , \\  
         & \\
       \frac{p^2 - p^2 \chi(p) - 3 p +  p \chi(p) + 3 \chi(p)}{p^3- 3 p^2 + 3 p + p \chi(p) - 3 \chi(p)} 
         \quad & \text{if} \quad p \mid N-2 .
    \end{cases}
\end{equation}

This gives
\[
  G_N(s) = L(s+1, \chi) \prod_p \left(  1 + \chi(p)p^{-s-1} K_N(p) \right) .
\]
We see that $G_N(s)$ has an analytic continuation to $Re \; (s) > -1$ and
\[
 G_N (s) \ll N^{\varepsilon} \, T^{1/6} \qquad \text{for} \qquad Re \; (s) \ge -\frac{1}{2} , \quad
    |Im \; (s)| \le T .
\]
Applying Perron's formula and proceeding as above we find
\begin{equation} \label{200}
  \sum_{d \le D} G_N(d) = G_N(0) + O \left( N^{-1/100} \right) 
   =
   \frac{\pi}{4} \prod_p \left( 1 + \frac{\chi(p)}{p} K_N(p) \right) + O \left( N^{-1/100} \right) .
\end{equation}

Using \eqref{187}, \eqref{189}, \eqref{196} and \eqref{200} we find
\begin{equation} \label{204}
 \Gamma(N) = \frac{1}{16} \mathfrak S_R(N) +  O \left( N^{\varepsilon - 1/100} \right) ,
\end{equation}
where $\mathfrak S_R(N)$ is defined by \eqref{50}.
We leave the calculations to the reader.

From \eqref{173.1}, \eqref{173.7}, \eqref{174} and \eqref{204}
we get
\begin{equation} \label{209}
  R_{1, 1}(N) 
    = 
  \frac{1}{32} N^2 \mathfrak S_R (N) + O (N^2 \mathcal L^{-1})  .
\end{equation}

\paragraph{The estimation of $R_{1, 2}(N)$.}

Using \eqref{52}, \eqref{56}, \eqref{171} and \eqref{173} we write 
\begin{equation} \label{210}
  R_{1, 2}(N) 
     = 
   \sum_{2 < p < N} \, (\log p) \, r_2(p - 1) \, \sum_{d \le D} \chi(d)
  I^{(2)}_{d, 1} (N-p)   = R_{1, 2}'(N) + R_{1, 2}^*(N) ,
\end{equation}
where
\begin{align}
 R_{1, 2}'(N) 
 & =
   \sum_{2 < p < N} \, (\log p) \, r_2(p - 1) \, \sum_{d \le D} \frac{\chi(d)}{\varphi(d)} \,
  \mathfrak S^{(2)}_{d, 1}(N-p) \, (N-p-1) ,
  \label{250} \\
  &   \notag \\
  R_{1, 2}^*(N)
  &   =
  \sum_{2 < p < N} \, (\log p) \, r_2(p - 1) \, \sum_{d \le D} \chi(d) \,
  \Delta^{(2)}_{d, 1}(N-p) .
  \label{251}
\end{align}

From \eqref{171}, \eqref{251} and Cauchy's inequality we find
\begin{align}
 |R_{1, 2}^*(N)| 
  &  \ll
  \mathcal L \sum_{2 < p < N}  \tau(p - 1) \, \sum_{d \le D} |
  \Delta^{(2)}_{d, 1}(N-p) | 
    \notag \\
    & \notag \\
  & \ll
    \mathcal L \sum_{n \le  N}  \tau(n) \, \sum_{d \le D} |
  \Delta^{(2)}_{d, 1}(n) |
    \notag \\
    & \notag \\
  & \ll
  \mathcal L \left( \sum_{n \le  N}  \, \sum_{d \le D} 
  \tau^2(n) |
  \Delta^{(2)}_{d, 1}(n) |\right)^{1/2} 
  \left( \sum_{n \le  N}   \, \sum_{d \le D} |
  \Delta^{(2)}_{d, 1}(n) |\right)^{1/2} 
  \notag \\
  & \notag \\
  &  = 
  \mathcal L \; U^{1/2}\; V^{1/2} ,
  \label{252}
\end{align}
say. We use the trivial bound  $\Delta^{(2)}_{d, 1}(n) \ll \mathcal L^2 N  d^{-1}$ and the well-known elementary inequa\-lity 
$\sum_{n \le x} \tau^2(n) \ll x \log^3 x$ and we find
\begin{equation} \label{252.5}
 U \ll N^2 \mathcal L^6 .
\end{equation}
To estimate $V$ we apply \eqref{150} and Lemma~\ref{L1} and we get
\begin{equation} \label{253}
 V \ll N^2 \mathcal L^{-10} .
\end{equation}
From \eqref{252}, \eqref{252.5} and \eqref{253} it follows that
\begin{equation} \label{254}
  R_{1, 2}^*(N)  \ll N^2 \mathcal L^{-1} .
\end{equation}

Consider now $R_{1, 2}'(N)$. Using \eqref{12}, \eqref{54}, \eqref{55} and \eqref{250}
we write it in the form
\[
  R_{1, 2}'(N) = c_0 \sum_{2 < p < N} (\log p) \, r_2(p-1) \, (N-p-1) \, \lambda(N-p)
  \sum_{\substack{ d \le D \\ (d, N-p-1) = 1 }} \frac{\chi(d)\, \lambda(d)}{\varphi(d) \, \lambda((d, N-p))} .
\]
It is not difficult to find an asymptotic formula for the sum over $d$. However such a formula is already established
in section 3.2 of \cite{Tolev1} and it implies that $\sum_d \ll \log \mathcal L$. Therefore, using also \eqref{12}, we find
\[
 R_{1, 2}'(N) \ll N \mathcal L (\log \mathcal L)^2 \sum_{p \le N } |r_2(p-1)| .
\]
It remains to apply \eqref{171} and Lemma~\ref{L4} and we get
\begin{equation} \label{256}
  R_{1, 2}'(N)  \ll N^2 \mathcal L^{-\theta_0} (\log \mathcal L)^7 .
\end{equation}
From \eqref{210}, \eqref{254} and \eqref{256} we obtain
\begin{equation} \label{258}
 R_{1, 2}(N)  \ll N^2 \mathcal L^{-\theta_0} (\log \mathcal L)^7 .
\end{equation}

\paragraph{The estimation of $R_{1, 3}(N)$.}

We use \eqref{60}, \eqref{171} and \eqref{173} to write
\begin{align}
  R_{1, 3}(N) 
    & = 
    \sum_{p_1 + p_2 + p_3 = N} (\log p_1 ) ( \log p_2 ) (\log p_3 )
     \sum_{\substack{d \mid p_1 - 1 \\ d \le D}} \chi(d) \; 
     \sum_{\substack{ m \mid p_2 - 1 \\ \frac{p_2 - 1}{ m } \ge N/D  }} \chi \left( \frac{p_2 - 1}{m} \right) 
    \notag \\
    & \notag \\
    & =
      \sum_{\substack{ d \le D \\ m \le D \\ 2 \mid m}} \chi(d) \sum_{j = \pm 1 } \chi(j) 
     I^{(3)}_{\langle d, 4 m \rangle , \langle 1 , 1 + j m \rangle}(N , \langle [1, N], J_{m} \rangle ) ,
  \notag 
\end{align}  
where $J_m = [1 + m N/D , N]$. 
From \eqref{75} we get
\begin{equation} \label{300}
  R_{1, 3}(N) = R_{1, 3}'(N) + R_{1, 3}^*(N) ,
\end{equation}
where
\begin{align}
 R_{1, 3}'(N)
   &  = 
  \sum_{\substack{ d \le D \\ m \le D \\ 2 \mid m }} 
    \frac{\chi(d) \; \Phi^{(3)} (N  , \langle [1, N], J_m \rangle )}{\varphi(d) \; \varphi(4 m)}  
  \sum_{j = \pm 1 }
   \chi(j) \;
   \mathfrak S^{(3)}_{\langle d , 4 m  \rangle , \langle 1 , 1 + j m \rangle}(N) ,
   \label{310} \\
R_{1, 3}^*(N)
   &  = 
  \sum_{\substack{ d \le D \\ m \le D \\ 2 \mid m }} 
    \chi(d) 
  \sum_{j = \pm 1 }
   \chi(j) \;
    \Delta^{(3)}_{\langle d , 4 m  \rangle , \langle 1 , 1 + j m \rangle} (N, \langle [1, N] , J_m \rangle) .
   \notag 
\end{align}

From \eqref{150} and Lemma~\ref{L2} we find
\begin{equation} \label{330}
  R_{1, 3}^*(N) \ll N^2 \mathcal L^{-1} .
\end{equation}
Consider $R_{1, 3}'(N)$. According to Lemma~\ref{L3} the expression
$\mathfrak S^{(3)}_{\langle d , 4 m  \rangle , \langle 1 , 1 + j m \rangle}(N) $
does not depend on $j$. Therefore from \eqref{310} it follows that
\begin{equation} \label{340}
  R_{1, 3}'(N) = 0 .
\end{equation}
Using \eqref{300}, \eqref{330} and \eqref{340} we obtain
\begin{equation} \label{480}
 R_{1, 3} (N) \ll
  N^2 \mathcal L^{ - 1 } .
\end{equation}

\paragraph{The estimation of $R_{2, 2}(N)$.}

Let $\mathcal P$ be the set of primes, specified in Lemma~\ref{L6} 
(with $\omega = B(10) + C(10)+ 1$). Using \eqref{150}, \eqref{171}, 
\eqref{173} and the inequality
$uv \le u^2 + v^2$ we get
\[
 R_{2, 2}(N) \ll \mathcal L^3 
   \sum_{\substack{ p_1 + p_2 + p_3 = N \\ p_2 \in \mathcal P }}
    \left| \sum_{\substack{ d \mid p_1 - 1 \\ D < d < N/D}} \chi(d) \right|^2
    = \mathcal L^3 
   \sum_{ p_1 < N }
    \left| \sum_{\substack{ d \mid p_1 - 1 \\ D < d < N/D}} \chi(d) \right|^2
    \sum_{\substack{ p_2 + p_3 = N - p_1 \\ p_2 \in \mathcal P }} 1 .
\]
We estimate the sum over $p_2, p_3$ using Lemma~\ref{L6} and we find
\[
 R_{2, 2}(N) \ll N \mathcal L^{1 - 2 \theta_0} (\log \mathcal L)^6
    \sum_{ p < N }
    \left| \sum_{\substack{ d \mid p - 1 \\ D < d < N/D}} \chi(d) \right|^2 .
\]
Then we apply Lemma~\ref{L5} and obtain
\begin{equation} \label{500}
 R_{2, 2} (N) \ll
  N^2 \mathcal L^{ - 2 \theta_0} (\log \mathcal L)^{13} .
\end{equation}

\paragraph{The estimation of $R_{2, 3}(N)$.}

Using \eqref{52}, \eqref{171} and \eqref{173} we write $R_{1, 2}(N)$ in the form
\begin{align} 
  R_{2, 3}(N) 
    &  = 
    \sum_{p_1 + p_2 + p_3 = N} (\log p_1) (\log p_2) (\log p_3) \;
      r_2(p_1 - 1) \sum_{\substack{m \mid p_2 - 1 \\ \frac{p_2 - 1}{m} \ge N/D}} 
      \chi \left(  \frac{p_2 - 1}{m}\right)      
      & \notag \\
     & =
     \sum_{2 < p < N}  (\log p) \; r_2 (p-1) 
     \sum_{\substack{m \le D \\ 2 \mid m}} \; \sum_{j = \pm 1} \chi(j) \;
       I^{(2)}_{  4m  ,  1 + j m   } (N-p , J_m) ,
      \notag
\end{align}      
where $J_m = [1 + m N / D, N]$.
Having in mind \eqref{56} we write
\begin{equation} \label{602}
  R_{2, 3}(N) = R'_{2, 3}(N) + R^*_{2, 3}(N) ,
\end{equation}
where
\begin{align}
 R_{2, 3}'(N) 
 & =
   \sum_{2 < p < N} \, (\log p) \, r_2(p - 1) \, 
       \sum_{\substack{m \le D \\ 2 \mid m}} \frac{\Phi^{(2)} (N-p \, , J_m)}{\varphi(4m)} \sum_{j = \pm 1} \chi(j) \,
      \mathfrak S^{(2)}_{4m, 1 + jm}(N-p) 
  \label{603} \\
  &   \notag \\
  R_{2, 3}^*(N)
  &  =
  \sum_{2 < p < N} \, (\log p) \, r_2(p - 1) \,  \sum_{\substack{m \le D \\ 2 \mid m}} \; \sum_{j = \pm 1}
    \chi(j) \;
  \Delta^{(2)}_{4m, 1 + jm}(N-p , J_m)
  \label{604}
\end{align}

Consider first $R_{2, 3}'(N)$. From Lemma~\ref{L3} we know that $\mathfrak S^{(2)}_{4m, 1 + jm}(N-p) $
does not depend on $j$. Hence using \eqref{603} we get
\begin{equation} \label{605}
  R'_{2, 3}(N) = 0 .
\end{equation}

Consider now $R_{2, 3}^*(N)$.
From \eqref{171}, \eqref{604} and Cauchy's inequality we find
\begin{align}
 R_{2, 3}^*(N)
  &  \ll
  \mathcal L \sum_{2 < p < N}  \tau(p - 1) \, \sum_{\substack{m \le D \\ 2 \mid m}} \sum_{j = \pm 1}
  | \Delta^{(2)}_{4m, 1 + jm}(N-p , J_m) | 
    \notag \\
    & \notag \\
  & \ll  
     \mathcal L \sum_{n \le  N}  \tau(n) \, \sum_{\substack{m \le D \\ 2 \mid m}} \sum_{j = \pm 1}
  | \Delta^{(2)}_{4m, 1 + jm}(n , J_m) | 
    \notag \\
    & \notag \\
  & \ll
   \mathcal L \; U_1^{1/2}\; V_1^{1/2} ,
  \label{606}
\end{align}
where
\[ 
 U_1  =
   \sum_{n \le  N}  
  \tau^2(n) \sum_{\substack{m \le D \\ 2 \mid m}} \; \sum_{j = \pm 1}
  \left| \Delta^{(2)}_{4m, 1 + jm}(n , J_m) \right| , \qquad
  V_1  =
   \sum_{n \le  N}   \, \sum_{\substack{m \le D \\ 2 \mid m}} \; \sum_{j = \pm 1}
  \left| \Delta^{(2)}_{4m, 1 + jm}(n , J_m) \right| .
\]  
We use the trivial bound  $\Delta^{(2)} \ll \mathcal L^2 N m^{-1}$ and the inequality 
$\sum_{n \le x} \tau^2(n) \ll x \log^3 x$ to find
\begin{equation} \label{610}
 U_1 \ll N^2 \mathcal L^6 .
\end{equation}
We estimate $V_1$ using \eqref{150} and Lemma~\ref{L1} and we get
\begin{equation} \label{612}
 V_1 \ll N^2 \mathcal L^{-10} .
\end{equation}
Using \eqref{606}, \eqref{610} and \eqref{612} we obtain
\begin{equation} \label{620}
  R_{2, 3}^*(N)  \ll N^2 \mathcal L^{-1} .
\end{equation}

Now taking into account \eqref{602}, \eqref{605} and \eqref{620} we find
\begin{equation} \label{630}
  R_{2, 3}(N)  \ll N^2 \mathcal L^{-1} .
\end{equation}

\paragraph{The estimation of $R_{3, 3}(N)$.}

We use \eqref{60}, \eqref{171} and \eqref{173} to write
\begin{align}
 R_{3, 3}(N) 
   & = 
    \sum_{p_1 + p_2 + p_3 = N} (\log p_1 ) ( \log p_2 ) (\log p_3 )
     \sum_{\substack{ m_1 \mid p_1 - 1 \\ \frac{p_1 - 1}{m_1} \ge N/D}} \chi \left( \frac{p_1 - 1 }{m_1}\right) 
     \sum_{\substack{ m_2 \mid p_2 - 1 \\ \frac{p_2 - 1}{m_2} \ge N/D}} \chi \left( \frac{p_2 - 1 }{m_2}\right) 
      \notag \\
     & \notag \\
   & = 
   \sum_{\substack{m_1, m_2 \le D \\ 2 \mid m_1 , \; 2 \mid m_2}} \; \sum_{\substack{j_1 = \pm 1 \\ j_2 = \pm 1}}
     \chi(j_1) \chi(j_2) \;
    I^{(3)}_{\langle 4m_1, 4 m_2 \rangle , \langle 1 + j_1 m_1 , 1 + j_2 m_2 \rangle}
    (N , {\bf J_m}) ,
      \notag 
\end{align}      
where ${\bf J_m} = \langle J_{m_1} , J_{m_2} \rangle $; $J_{m_{\nu}} = [ 1+ m_{\nu} N/ D, N]$, $\nu = 1, 2 $.
We write
\begin{equation} \label{700}
  R_{3, 3}(N)  = R_{3, 3}'(N) +  R_{3, 3}^*(N)  ,
\end{equation}
where
\begin{align}
 R_{3, 3}'(N)
  &  = 
  \sum_{\substack{m_1, \; m_2 \le D \\ 2 \mid m_1 , \; 2 \mid m_2}} 
  \frac{\Phi^{(3)} (N , {\bf J_m})}{\varphi(4 m_1) \varphi(4 m_2)}
     \sum_{\substack{j_1 = \pm 1 \\ j_2 = \pm 1}}
   \chi(j_1) \chi(j_2) \;
 \mathfrak S^{(3)}_{\langle 4m_1, 4 m_2 \rangle \, , \, \langle 1 + j_1 m_1 , 1 + j_2 m_2 \rangle}(N) ,
    \label{710} \\
    & \notag \\
R_{3, 3}^*(N)
  &  = 
  \sum_{\substack{m_1, \; m_2 \le D \\ 2 \mid m_1 , \; 2 \mid m_2}} 
    \;
     \sum_{\substack{j_1 = \pm 1 \\ j_2 = \pm 1}}
   \chi(j_1) \chi(j_2) \;
 \Delta^{(3)}_{\langle 4m_1, 4 m_2 \rangle \, , \, \langle 1 + j_1 m_1 , 1 + j_2 m_2 \rangle}(N, {\bf J_m}) .
    \notag
\end{align}
According to Lemma~\ref{L3}, the expression $\mathfrak S^{(3)}$ in \eqref{710} does not depend of 
$j_2$ and therefore
\begin{equation} \label{760}
    R_{3, 3}'(N) =0 .
\end{equation}
On the other hand, using \eqref{150} and Lemma~\ref{L2} we get
\begin{equation} \label{750}
    R_{3, 3}^*(N) \ll N^2 \mathcal L^{-1} .
\end{equation}

From \eqref{700}, \eqref{750} and \eqref{760} it follows that
\begin{equation} \label{1000}
 R_{3, 3} (N) \ll
  N^2 \mathcal L^{ - 1 } .
\end{equation}

\paragraph{The end of the proof.}
The asymptotic formula \eqref{40} is a consequence of 
\eqref{172}, \eqref{173}, \eqref{209}, \eqref{258}, \eqref{480}, \eqref{500}, \eqref{630} and \eqref{1000}. 
The theorem is proved.
$\square$

\bigskip
\bigskip

\vbox{
\hbox{Faculty of Mathematics and Informatics}
\hbox{Sofia University ``St. Kl. Ohridsky''}
\hbox{5 J.Bourchier, 1164 Sofia, Bulgaria}
\hbox{ }
\hbox{Email: dtolev@fmi.uni-sofia.bg}}

\end{document}